\documentclass[a4paper,twoside,12pt]{amsart}

\usepackage[utf8]{inputenc}
\usepackage[T1]{fontenc}
\usepackage[english]{babel}
\usepackage{tgpagella}
\usepackage{csquotes}
\usepackage{fancyhdr}
\usepackage{multirow}
\usepackage{amsthm, amsfonts, amssymb, amsmath, latexsym, enumerate, array}
\usepackage{amscd} 
\usepackage[all]{xy}
\usepackage{pstricks,graphicx}
\usepackage{stmaryrd}
\usepackage{pgf,tikz}
\usetikzlibrary{arrows}
\usepackage{url}
\usepackage{braket}
\usepackage{hyperref,mdwlist,mathrsfs}

\usepackage{color}
\usepackage{paralist}

\newcommand{\PP}{\ensuremath{\mathbb{P}}}
\newcommand{\HF}{\ensuremath{\mathrm{HF}}}

\newtheorem{theorem}{Theorem}[section]
\newtheorem{lemma}[theorem]{Lemma}

\newtheorem{proposition}[theorem]{Proposition}

\newtheorem*{proposition*}{Proposition}
\theoremstyle{definition}
\newtheorem{definition}[theorem]{Definition}
\newtheorem{example}[theorem]{Example}
\newtheorem{remark}[theorem]{Remark}
\newtheorem{problem}{Problem}

\DeclareMathOperator{\Ann}{Ann}
\DeclareMathOperator{\rank}{rank}
\DeclareMathOperator{\hess}{hess}

\DeclareMathOperator{\Tor}{Tor}
\DeclareMathOperator{\reg}{reg}

\title[Perazzo $n$-folds and the weak Lefschetz property]{Perazzo $n$-folds and the weak Lefschetz property}

 \author[]{Emilia Mezzetti} 
 \address{\hspace{-15pt} Dipartimento di Matematica, Informatica e  Geoscienze, Universit\`a di
Trieste, Via A.Valerio 12/1, 34127 Trieste, Italy}
  \email{mezzette@units.it, ORCID 0000-0001-5300-9779}
 
\author[]{Rosa M.\ Mir\'o-Roig} 
  \address{\hspace{-15pt} Facultat de
  Matem\`atiques i Inform\`atica, Universitat de Barcelona, Gran Via des les
  Corts Catalanes 585, 08007 Barcelona, Spain} \email{miro@ub.edu, ORCID 0000-0003-1375-6547}

\thanks{\hspace{-15pt}  Mezzetti is a member of INdAM - GNSAGA, Mir\'o-Roig was partially supported by the grant PID2020-113674GB-I00.
}

\begin{document}

\begin{abstract} 
In this paper, we determine the maximum $h_{max}$ and the minimum $h_{min}$ of the Hilbert vectors of Perazzo algebras $A_F$, where $F$ is a Perazzo polynomial of degree $d$ in $n+m+1$ variables. These algebras always fail the Strong Lefschetz Property. We determine the integers $n,m,d$ such that $h_{max}$ (resp. $h_{min}$) is unimodal, and we prove that $A_F$ always fails the Weak Lefschetz Property if its Hilbert vector is maximum, while it satisfies the Weak Lefschetz Property if it is minimum, unimodal, and satisfies an additional mild condition. We determine the minimal free resolution of Perazzo algebras associated to Perazzo threefolds in $\PP^4$ with minimum Hilbert vectors. Finally we pose some open problems in this context.

\end{abstract}
\maketitle
{\it Dedicated to Enrique Arrondo on the occasion of his $60^{th}$ birthday}

\section{Introduction}

A Perazzo form of degree $d$ is by definition (see \cite{Go}) a homogeneous polynomial $F\in K[X_0, \ldots,X_n,U_1,\ldots,U_m]$
\[
F=X_0p_0+X_1p_1+\dots+X_np_n+G,
\]
with $n\geq m\geq 2$, $p_0, \ldots, p_n\in K[U_1, \ldots, U_m]_{d-1}$,  $G\in K[U_1, \ldots, U_m]_d$, where $p_0, \ldots, p_n$ are algebraically dependent but linearly independent.

The interest of Perazzo forms comes form the fact that their sets of zeros are a class of examples of hypersurfaces with vanishing hessian which are not cones. The problem of characterizing the hypersurfaces not cones with vanishing hessian is an important open problem in projective algebraic geometry since the classical works of Hesse and of Gordan--Noether (\cite{He1}, \cite{He2}, \cite{GN}). Hesse believed that the hessian of a homogenous polynomial  vanishes if and only if its variety of zeros  is
a cone, as is indeed the case when the degree of the form is $2$. However, P. Gordan and M. Noether 
proved that while Hesse's claim is true for forms in at most $4$ variables, it is false for $5$ or more variables and any
degree $\geq 3$. They also gave a complete description of the hypersurfaces in $\mathbb P^4$, non cones, having vanishing hessian: their equations are elements of $K[U_1, U_2][\Delta]$, where $\Delta$ is a Perazzo polynomial of the form $X_0p_0+X_1p_1+X_2p_2$. Since then, many efforts have been made to find a characterization for arbitrary degree and number of variables .  It turns out that all the counterexamples known so far can be built starting from Perazzo forms (see \cite[ Appendix A]{Go}, \cite[Chapter 7]{R}). The original paper of Perazzo \cite{Pe} deals with the case of cubic hypersurfaces, it was revisited in \cite{GR}.

The study of hessians of homogeneous polynomials has gained new attention because of its connection to Lefschetz properties for graded Artinian Gorenstein algebras. Recall
that a standard graded Artinian algebra $A$ has the weak Lefschetz property (WLP) if
multiplication by a generic linear form $\ell$ has maximal rank in each degree. Similarly
$A$ has the strong Lefschetz property (SLP) if multiplication by $\ell^s$ has maximal rank
in each degree for every positive integer $s$. The study of the Lefschetz properties for graded Artinian algebras originates from the Hard Lefschetz Theorem, which implies that the cohomology ring $A$ of any smooth complex projective variety has the SLP;
moreover $A$ is a graded Artinian Gorenstein algebra. 

Although Lefschetz properties have been the subject of intense research in recent years, many natural problems are still open and the general picture is far from being understood. In the Gorenstein case every standard graded Artinian algebra can be written as $A_F$, the quotient of a ring of differential operators by the annihilator of a homogeneous polynomial $F$,  called its Macaulay dual generator (see Section \ref{subsection Hilbert functions} for details). 
Due to work of Watanabe and Maeno -- Watanabe (\cite{W}, \cite{MW}), a non-trivial characterization of Artinian Gorenstein algebras failing the SLP, in terms of the Macaulay dual generator, has been found.
Indeed,  $A_F$ fails the SLP if and only if one of the
non-trivial higher hessians of $F$ vanishes. This result has been generalised to the WLP using the so called
mixed hessians (see \cite{GZ}).

It follows that the Artinian Gorenstein algebras $A_F$ associated to Perazzo forms fail the SLP. 
It is therefore natural to pose the question if these algebras satisfy or fail the WLP. This question has been considered in some recent articles (\cite{FMMR22}, \cite{MRP}, \cite{A}), where the case of Perazzo forms with $m=2$ has been completely solved.

Before summarizing the results of those articles, we recall a few basic facts about the Hilbert functions of graded Artinian Gorenstein algebras. If $A$ is such an algebra of socle degree $d$, then its Hilbert function is captured by its $h$-vector $(h_0, h_1, \dots, h_d)$, where $h_i=\dim_K A_i$. Since $A$ is a Poincar\' e duality algebra, the $h$-vector results to be symmetric, i.\,e.\ $h_i=h_{d-i}$. On the set of $h$-vectors of the same length there is the natural componentwise partial order: given  $h = (h_0, h_1, \dots, h_d)$ and 
$h'=(h'_0, h'_1,\dots,h'_d)$, we say that $h \leq h'$ if
$h_i \leq h'_i$
 for every $i$, $0 \leq i \leq d$.

In the quoted articles the following facts are proved. Let $A_F$ be the Artinian Gorenstein algebra associated to a Perazzo form $F$ with $n\geq m=2$, $d\geq n+1$. Let $(h_0, h_1,\ldots, h_d)$ be its $h$-vector. Then:
\begin{enumerate}
    \item the Hilbert function of $A_F$ is unimodal, i.\,e.\ $h_0\leq h_1 \leq \dots \leq h_k \geq h_{k+1}\geq \dots \geq h_d$ for some $k$;
    \item $h_i\leq d+2$ for any index $i$ and $A_F$ has the WLP if and only if   \  $\sharp\{i |h_i=d+2\}\leq 1$;
    \item the $h$-vectors of the algebras $A_F$ for fixed $n,m,d$ have a maximum and a minimum, that are completely described. In particular if the $h$-vector is maximum $A_F$ fails the WLP, while if the $h$-vector is minimum $A_F$ has the WLP if and only if $d\geq 2n$;
    \item the Perazzo forms in $5$ variables such that the $h$-vector is minimum admit a precise description in terms of the position of the $2$-plane generated by $p_0, p_1, p_2$ in $\PP (K[U_1, U_2]_{d-1})$ with respect to the rational normal curve.
\end{enumerate}

In this article we consider Perazzo forms in any number of variables with the aim to extend the results obtained for $m=2$. Moreover we tackle the problem of describing the minimal free resolutions of Perazzo algebras. We are able to prove that, for fixed $m,n,d$, the $h$-vectors of Perazzo algebras have a maximum and a minimum as in the case $m=2$, and we describe them in Proposition \ref{partial} and Theorem \ref{minimal} respectively. But, diversely from the case $m=2$, if $m\geq 3$ these $h$-vectors are not always unimodal. In Theorems \ref{non unimodal},  \ref{non unimodal2} and \ref{unimodality2} we characterize the integers $n\geq m\geq 3$ such that $h_{max}$ (resp. $h_{min}$) is not unimodal. We note that $h_{max}$ is never unimodal for $d$ large enough.

Regarding the WLP,  we find that Perazzo algebras with maximal $h$-vector never have the WLP, while those with minimal $h$-vector have the WLP provided that $h_{min}$ is unimodal and an additional mild condition is satisfied. 
The problem of characterizing when WLP holds for intermediate $h$-vectors remains open.

As for our second aim, we are able to compute in Theorem \ref{main} the minimal free resolution for a class of Perazzo algebras, those in $5$ variables with minimal $h$-vector. The proof is by induction on the degree $d$, the base of the induction being possible because of the explicit description of the algebras $A_F$ with minimal $h$-vector.

Many questions remain open, and we devote the last section of this article to list a few open problems that we think deserve to be considered.
\medskip

The paper is organized as follows. We start by reviewing in Section \ref{background} definitions and basic results concerning Artinian Gorenstein algebras associated to  Perazzo hypersurfaces, minimal free resolutions and Lefschetz properties. In Section \ref{maximal}, we determine the maximal Hilbert function once the integers $n,m,d$ are fixed. We study when this function is unimodal and we prove that Perazzo algebras with maximal Hilbert function do not have the WLP. In Section \ref{section minimal}, similarly, we determine the minimal Hilbert function, study its unimodality and prove that Perazzo algebras with this Hilbert function, and satisfying an additional mild condition that implies the unimodality of the $h$-vector, have the WLP. We also characterize the integers $n,m,d$ such that the maximum and the minimum Hilbert function coincide. In Section \ref{MFRPerazzo} we compute the
minimal free resolution for the Perazzo algebras in $5$ variables with minimal $h$-vector.
Finally, in Section \ref{final} we pose some relevant open problems in this circle of ideas.

\vskip 4mm
\textbf{Acknowledgement.} The Authors thank Luca Fiorindo, Pedro Macias Marques and Lisa Nicklasson  for interesting discussions, and the anonymous referee for useful comments. 

\section{Background}\label{background}
In this section we fix notations, we recall the basic facts on Hilbert functions, Lefschetz properties, minimal free resolutions as well as on  Perazzo hypersurfaces needed later on.

\subsection{Hilbert functions}\label{subsection Hilbert functions}
Throughout this paper  $K$ will be an algebraically closed field of characteristic zero.
Given a standard graded Artinian $K$-algebra $A=R/I$ where $R=K[x_0,x_1,\dots,x_N]$ and $I$ is a homogeneous ideal of $R$,
we denote by $\HF_A\colon \mathbb{Z} \longrightarrow \mathbb{Z}$ with $\HF_A(j)=\dim _KA_j=\dim _K[R/I]_j$
its Hilbert function. Since $A$ is Artinian, its Hilbert function is
captured in its \emph{$h$-vector} $h=(h_0,h_1,\dots ,h_d)$ where $h_i=\HF_A(i)>0$ and $d$ is the last index with this property. The integer $d$ is called the \emph{socle degree of} $A$.  We will use the terms ``Hilbert function'' and ``$h$-vector'' interchangeably along the paper.

We recall the construction of the Artinian Gorenstein algebra $A_F$ with Macaulay dual generator a given form $F\in S=K[X_0,\dots,X_N]$; we denote by $R=K[x_0,\ldots,x_N]$ the ring of differential operators acting on the polynomial ring $S$, i.\,e.\ $x_i=\frac{\partial}{\partial X_i}$. Therefore $R$ acts on $S$ by differentiation. Given polynomials $p\in R$ and $G\in S$ we will denote by ${p\circ G}$ the differential operator $p$ applied to $G$. We define
\[
\Ann_RF:=\{p\in R \mid p\circ F=0\}\subset R,
\]
and $A_F=R/\Ann_R F$: it is a standard graded Artinian Gorenstein $K$-algebra and $F$ is called its Macaulay dual generator. We remark that every standard graded Artinian Gorenstein $K$-algebra is of the form $A_F$ for some form $F$, in view of the \lq\lq Macaulay double annihilator Theorem'' (see for instance \cite[Lemma 2.12]{IK}). We may abbreviate and write ${\Ann F}$ when the ring $R$ is understood.

As an important key tool to determine the unimodality of the Hilbert function of a Perazzo algebra or the minimal free resolution of  Artinian Gorenstein algebras associated to Perazzo threefolds in $\PP^4$, we state the following:

\begin{proposition}
\label{exactsequence}
Let $A_F$ be an Artinian Gorenstein graded $K$-algebra and set $I=\Ann F$. Then for every linear form $\ell\in A_1$ the sequence
\begin{equation}
0\longrightarrow \frac{R}{(I\colon \ell)}(-1)\longrightarrow A_F=\frac{R}{I}\longrightarrow \frac{R}{(I,\ell)}\longrightarrow 0
\end{equation}
is exact. Moreover $\frac{R}{(I\colon\ell)}$ is an Artinian Gorenstein graded algebra with $\ell\circ F$ as Macaulay dual generator.
\end{proposition}
\begin{proof} 
We get the result cutting the exact sequence 
\[
0\longrightarrow \frac{(I\colon\ell)}{I}(-1)
\longrightarrow \frac{R}{I}(-1)\xrightarrow{\,\,\,\times\ell\,\,\,}\frac{R}{I}\longrightarrow \frac{R}{(I,\ell)}\longrightarrow 0
\]
into two short exact sequences. The second fact is a straightforward computation.
\end{proof}

\subsection{Minimal free resolutions.} Let $A=R/I$ be an Artinian graded $K$-algebra. It is well known that it has a minimal graded free $R$-resolution of the following type:
$$
0\longrightarrow F_{n+1} 
\longrightarrow F_n \longrightarrow \cdots  \longrightarrow F_i \longrightarrow \cdots \longrightarrow F_1\longrightarrow  R\longrightarrow A\longrightarrow 0 
$$
where
$$ F_i=\oplus _jR(-j)^{\beta_{ij}^R(A)}
$$
and the graded Betti numbers $\beta_{ij}^R(A)$ of $A$ over $R$ are defined as usual as  the integers 
$$\beta_{ij}^R(A)=\dim_K [\Tor^R_i(A,K)]_j.$$ These homological invariants are our main focus and indeed our goal in Section \ref{MFRPerazzo} is to determine the graded Betti numbers $\beta _{ij}^R(A_F)$ of an Artinian Gorenstein algebra $A_F$ associated to a Perazzo 3-fold $F$ in $\PP^4$ with termwise minimal Hilbert function. It is important to point out that the graded Betti numbers of an Artinian graded algebra $A$ determine its  graded Poincar\'{e} series. In fact, the graded Poincar\'{e} series of $A$ over $R$ is the generating function $$P^R_A(t,s)=\sum_{i,j} \beta_{ij}^R(A) t^is^j.$$ If $R$ is regular, then the Poincar\'{e} series is in fact a polynomial.

\subsection{Perazzo hypersurfaces}
\label{sec:prel_perazzo}

The simplest counterexample to Hesse's claim, i.\,e.\ a form with vanishing  hessian which does not define a cone, is $XU^2+YUV+ZV^2$. This example was extended to a class of cubic counterexamples in all dimensions  by Perazzo in \cite{Pe}. 

\begin{definition}
\label{perazzo} 
A \emph{Perazzo hypersurface}  $X=V(F) \subset \PP^N$ is the hypersurface defined by a \emph{Perazzo form}
\[
F=X_0p_0+X_1p_1+\cdots +X_np_n+G \in K[X_0,\ldots ,X_n,U_1\ldots ,U_m]_d
\]
where $n,m\ge 2$, $N=n+m$, $p_i\in K[U_1,\ldots ,U_m]_{d-1}$ are algebraically dependent but linearly independent, and $G\in K[U_1,\ldots ,U_m]_{d}$.

The Artinian Gorenstein algebra $A_F$ associated to a Perazzo polynomial will be called Perazzo algebra.
\end{definition}

The fact that the $p_i$'s are algebraically dependent implies $\hess_F=0$, while the linear independence assures that $V(F)$ is not a cone. 

We note that, to allow the linear independence of $p_0, \ldots, p_n$, we must assume 
\begin{equation}\label{inequality}
n+1\leq \binom{d+m-2}{m-1}.
\end{equation}
If equality holds in (\ref{inequality}) $p_0, \ldots, p_n$ form a basis of $K[U_1,\ldots,U_m]_{d-1}$; in this case $A_F$ is called a {\it full Perazzo algebra}. Full Perazzo algebras were studied in \cite{CGIZ} and  \cite{BGIZ}.
On the other hand, to guarantee the algebraic dependence for a general choice of $p_0, \ldots, p_n$ we make the assumption that $n\geq m$.

\vskip 4mm
The following lemma plays a key role in the induction step used in the proof of our main results (Theorem \ref{mainthm0} and Theorem \ref{main}).

\begin{lemma}
\label{lemmapartials2}
Let 
$F=X_0p_0+X_1p_1+\dots+X_np_n+G$ be a Perazzo form of degree $d$  and let $A_F$ be the associated Artinian Gorenstein algebra. Assume $n+1\leq \binom{d+m-3}{m-1}$. Then, for a general linear form $\ell\in A_F$, the polynomial $\ell\circ F$ defines a Perazzo form of degree $d-1$.
\end{lemma}
\begin{proof}
We can write $\ell=a_0X_0+a_1X_1+ \dots a_nX_n+b_1U_1+\dots +b_mU_m$ for some coefficients $a_i,b_j\in K$ not all zero. Then we can exhibit the action of $\ell$ on $F$ as
\begin{equation*}
\begin{split}
\ell\circ F=
X_0\Tilde{p_0}+ \dots + X_n\Tilde{p_n}+\left(a_0p_0+\dots+a_np_n+b_1\frac{\partial G}{\partial U_1}+\dots +b_m\frac{\partial G}{\partial U_m}\right)
\end{split}
\end{equation*}
with 
\begin{gather*}
\Tilde{p_0}=b_1\frac{\partial p_0}{\partial U_1}+\dots +b_m\frac{\partial p_0}{\partial U_m},\ldots, \Tilde{p_n}=b_1\frac{\partial p_n}{\partial U_1}+\dots +b_m\frac{\partial p_n}{\partial U_m}. 
\end{gather*}

The form $\ell\circ F$ has degree $d-1$. It remains to prove that the polynomials $\Tilde{p_0},\ldots,\Tilde{p_n}$ are linearly independent, for a general choice of $\ell$.  

Let ${\ell'=b_1U_1+\dots +b_mU_m}$ and consider the map 
\[\phi_b: {K[U_1, \ldots, U_m]_{d-1}\to K[U_1, \ldots, U_m]_{d-2}},\] given by ${p\mapsto \ell'\circ p}$. 
The fact that the  $(d-1)-$th powers  of linear forms span $K[U_1, \ldots, U_m]_{d-1}$ (\cite[Corollary 3.2]{I}), and similarly for $K[U_1, \ldots, U_m]_{d-2}$, implies that $\phi_b$ is surjective. 
Therefore its kernel 
has dimension  $\binom{m+d-3}{m-2}$. Now let ${W=\langle p_0,\ldots,p_n\rangle}$. Since the $p_i$ are linearly independent, $W$ has dimension $n+1$. It does not fill the whole space $K[U_1, \ldots,U_m]_{d-1}$, because $n+1<\binom{m+d-2}{m-1}$. Using again that the  $(d-1)-$th powers  of linear forms span $K[U_1, \ldots, U_m]_{d-1}$  and because of the assumption $n+1\leq \binom{m+d-3}{m-1}$, we deduce that there is an open dense set $B$ in $K^m$ such that for any ${(b_1, \ldots, b_m)\in B}$, the vector space $\ker \phi_b$ misses $W$. So, for any such $m$-tuple $(b_1,\ldots, b_m)$, the map 
\[
{\ell':W\to\ell'\circ W}
\]
is an isomorphism, and therefore $\Tilde{p_0},\ldots,\Tilde{p_n}$ are linearly independent.
\end{proof}

\subsection{Lefschetz properties}\label{subsection Lefschetz}

\begin{definition}\label{def:Lefschetz}
Let $A=R/I=\bigoplus_{i=0}^d\, A_i$ be a graded Artinian {$K$-algebra}. We say that $A$ has the \emph{weak Lefschetz property} (WLP, for short) if there is a linear form $\ell \in A_1$ such that, for all integers $i\ge0$, the multiplication map
\[
\times \ell\colon A_{i}  \longrightarrow  A_{i+1}
\]
has maximal rank, i.\,e.\ it is injective or surjective. In this case, the linear form $\ell$ is called a \emph{weak Lefschetz element} of $A$. We say that $A$ fails the WLP in degree $j$ if for a general form $\ell\in A_1$, the map $\times \ell\colon A_{j-1}  \longrightarrow  A_{j}$ does not have maximal rank. 

We say that $A$ has the \emph{strong Lefschetz property} (SLP, for short) if there is a linear form $\ell\in A_1$ such that, for all
integers $i\ge0$ and $k\ge 1$, the multiplication map
\[
\times \ell^k\colon A_{i}  \longrightarrow  A_{i+k}
\]
has maximal rank.  Such an element $\ell$ is called a \emph{strong Lefschetz element} of $A$.

\end{definition} 

It is easy to prove that the $h$-vector $(h_0, h_1, \ldots, h_d)$ of any graded Artinian $K$-algebra having the SLP or the WLP is \emph{unimodal}, i.e. there exists an index $k$ such that  $h_0 \leq h_1 \leq \dots \leq h_k \geq h_{k+1} \geq \dots \geq h_d$.

Let $A_F$ be an Artinian Gorenstein algebra associated to a Perazzo hypersurface of degree $d \ge 5$ in $\PP^4$. Recall that by \cite[Theorem 4.3]{FMMR22} the algebra $A_F$ has the weak Lefschetz property if the Hilbert function of $A_F$ is the termwise minimal one, namely $(1,5,6, \ldots, 6,5,1)$. Furthermore, in \cite[Proposition 3.7]{FMMR22} and \cite[Theorem 4.1]{FMMR22} it is proved that the maximal possible Hilbert function is
\begin{equation}\label{maxHF}
h_i=\left\{ 
\begin{aligned}
    &4i+1 \ \text{for } 1\le i\le\tfrac{d+1}{4}\\
    &d+2 \ \text{for } \tfrac{d+1}{4}<i\le \tfrac{d}{2}\\
    &\text{symmetry}
\end{aligned}
\right.
\end{equation}
and that any algebra $A_F$ with Hilbert function  as in (\ref{maxHF}) fails the WLP. 
As a complete  classification of Artinian Gorenstein algebras associated to Perazzo hypersurfaces of degree $d \ge 5$ in $\PP^4$ with the weak Lefschetz property we have the following result.
\begin{theorem}\label{thm:WLP} 
Let $A_F$ be an Artinian Gorenstein algebra associated to a Perazzo hypersurface $V(F)\subset \PP^4$ of degree $d\ge 5$. Let $(h_0,h_1,\ldots ,h_d)$ be its $h$-vector. The algebra $A_F$ has the WLP if and only if $\# \{i \mid h_i=d+2 \} \le 1$.
\end{theorem}
\begin{proof}
    See \cite[Theorem 3.11]{A}.
\end{proof}


\section{Maximal Hilbert function of a Perazzo algebra}\label{maximal}

In this section we determine the maximal $h$-vector $h_{max}=h_{max}(A_F)$ of a Perazzo algebra $A_F$ for any given $m, n, d$, extending the results obtained in the case $n=m=2$ in \cite{FMMR22} and \cite{A}; and $m=2$ and $n\ge 2$ in \cite{MRP}. Let $d$ be the degree of the Macaulay dual generator $F$ of $A_F$. We will see that, differently from the case $m=2$, for $m\geq 3$ and $d\gg 0$ $h_{max}$  is not  always unimodal.

Let $F$ be a Perazzo form as in Definition \ref{perazzo}:
\[
F=X_0p_0+X_1p_1+\cdots +X_np_n+G \in S_d=K[X_0,\ldots ,X_n,U_1\ldots ,U_m]_d,
\]
with $p_0,\ldots, p_n$ algebraically dependent but linearly independent.

\medskip

We will use the following notations: for $i=0,\ldots, n$ 
\begin{equation}\label{notation}
    p_i=\sum_{|\lambda |=d-1}\binom{d-1}{\lambda}p^i_{\lambda}U^\lambda
\end{equation}
where $\lambda=(\lambda_1, \ldots, \lambda_m)$ is a multi-index, $|\lambda|=\lambda_1+\dots+\lambda_m$, $\binom{d-1}{\lambda}= \frac{(d-1)!}{\lambda_1!\dots\lambda_m!}$ is the multinomial coefficient, and $U^\lambda=U_1^{\lambda_1}\dots U_m^{\lambda_m}$. 

Then, for any multi-index $\gamma=(\gamma_1, \ldots, \gamma_m)$ such that $|\gamma|\leq d-1$, the  partial derivative $\frac{\partial^{|\gamma|} p_i}{\partial U^\gamma}$ of $p_i$ with respect to $\gamma$ is equal to
\begin{equation*}
    (d-1)(d-2)\dots(d-1-|\gamma|)\sum_{|\mu|=d-1-|\gamma|}\binom{d-1-|\gamma|}{\mu}p^i_{\mu+\gamma}U^\mu.
\end{equation*}

      Similarly we put 
      $G=\sum_{|\lambda|=d}\binom{d}{\lambda}G_\lambda U^\lambda$.
      
      This notation will be useful to compute the $h$-vector of $A_F$, which is equivalent to computing  the dimension of $\Ann_R(F)_i$ for $i=0,\ldots, [\frac{d}{2}]$.

      \begin{proposition}\label{rank}    
          Let $i\leq [\frac{d}{2}]$ be an integer number. Let $h=(h_0,\ldots, h_d)$ be the $h$-vector of the Perazzo algebra $A_F$. Then $h_i$ is equal to the rank of the matrix containing in the columns the coefficients of the partial derivatives of $F$ of order $i$.  
      \end{proposition}
      \begin{proof}
      We put $R=K[x_0,\ldots,x_nu_1,\ldots, u_m]$. We observe that $h_1=n+m+1$. From now on we assume $i\ge 2$. Since $h_i=\dim (A_F)_i=\dim [R/\Ann_R(F)]_i$, we need to determine the polynomials $\phi$ of $R_i$ such that $\phi\circ F=0$.   Being $F$ linear in $x_0,\ldots,x_n$, all polynomials $\phi\in  R_i$ of degree at least $2$ in $x_0,\ldots, x_n$ clearly belong to $\Ann_R(F)_i.$ So assume that $\phi$ has degree $\leq 1$ in $x_0,\ldots, x_n$. We can write
      \[
      \phi=x_0\sum_{|\mu|=i-1}\alpha^0_\mu u^\mu+x_1\sum_{|\mu|=i-1}\alpha^1_\mu u^\mu+\dots+\sum_{|\nu|=i} \beta_\nu u^\nu.
      \]
      Imposing $\phi\circ F=0$ we get:
      \begin{equation}\label{conditions}
      \sum_{|\mu|=i-1}\alpha^0_\mu \frac{\partial^{|\mu|}p_0}{\partial U^{\mu}}+
      \sum_{|\mu|=i-1}\alpha^1_\mu \frac{\partial^{|\mu|}p_1}{\partial U^{\mu}}+\dots+
      X_0\sum_{|\nu|=i}\beta_\nu \frac{\partial^{|\nu|}p_0}{\partial U^\mu}+\dots +\sum_{|\nu|=i}\beta_\nu \frac{\partial^{|\nu|}G}{\partial U^\nu}=0.
      \end{equation}
      This gives rise to a homogeneous linear system of equations in the unknowns $\alpha_\mu^0, \alpha_\mu^1,\ldots, \beta_\nu$, with $|\mu|=i-1, |\nu|=i.$ The equations of the system are obtained equaling to zero the coefficients of the monomials of degree $d-i$ in $U_1, \ldots, U_m$ and in $X_0, \ldots, X_n, U_1, \ldots, U_m$. Being char $K= 0$, in view of notation (\ref{notation}), it follows that, up to a non-zero constant, the coefficients of the unknowns are precisely the coefficients of the partial derivatives of $F$ of order $i$. The theorem is proved.
      \end{proof}
\begin{proposition}\label{partial}
The maximal $h$-vector of a Perazzo algebra $A_F$, for fixed $m,n,d$ satisfying (\ref{inequality}), is $h_{max}=(h_0, \ldots, h_d)$ with \begin{equation*}
    h_i=\min\{\alpha_i+\beta_i, \alpha_i+\gamma_i\}
    \end{equation*}
      for any $0\leq i\leq [\frac{d}{2}]$, where
\begin{equation*}
 \alpha_i=\binom{m+i-1}{m-1}, \ \beta_i=\binom{d+m-i-1}{m-1} \text{ and } \gamma_i=(n+1)\binom{m+i-2}{m-1}.   
\end{equation*}

\end{proposition}
\begin{proof}
  It follows from  a result of Iarrobino, that we recall in Lemma \ref{Iarrobino}. Let $F$ be a Perazzo polynomial as in Definition \ref{perazzo} and assume that $p_0, \ldots, p_n, G$ are {\it general}. Let us compute the $h$-vector of $A_F$. In view of Proposition \ref{rank}, for any $i$, $1\leq i \leq [\frac{d}{2}]$, we have
  \begin{multline*}    
      h_i=\dim (A_F)_i=\\
      =\dim\langle\frac{\partial^iF}{\partial X_0^{i_0}\dots \partial X_n^{i_n}\partial U_1^{i_{n+1}}\dots \partial U_m^{i_{n+m}}} | \ i_0+\dots +i_{n+m}=i\rangle 
  \end{multline*}
  where $i_j\geq 0$ for $j=0,\ldots, n+m.$
 Being $p_0, \ldots, p_n$ general, this is equal to $\dim A+\dim B$ where
 \begin{equation*}
     A=\braket{ \frac{\partial^{i-1}p_0}{\partial U_1^{j_1}\dots \partial U_m^{j_m}}, \ldots, \frac{\partial^{i-1}p_n}{\partial U_1^{j_1}\dots \partial U_m^{j_m}} \mid j_1+\dots+j_m=i-1, j_r\geq 0 },
     \end{equation*}
\begin{multline*}
    B=\langle \frac{\partial^iF}{\partial U_1^{j_1}\dots \partial U_m^{j_m}}\mid j_1+\dots +j_m=i, j_r\geq 0 \rangle=\\
    =\langle \sum_{j=0}^n X_j \frac{\partial^ip_j}{\partial U_1^{j_1}\dots \partial U_m^{j_m}}+\frac{\partial^i G}{\partial U_1^{j_1}\dots \partial U_m^{j_m}}\mid j_1+\dots +j_m=i, j_r\geq 0\rangle.
\end{multline*}
    From Lemma \ref{Iarrobino} we get
    \begin{equation*}
        \dim A=\min\Set{ (n+1)\binom{m+i-2}{m-1},
         \binom{d+m-i-1}{m-1}}=\min\{\beta_i, \gamma_i\}, 
         \end{equation*}
        \begin{equation*}
        \dim B=\binom{m+i-1}{m-1}=\alpha_i,
    \end{equation*}
    which proves the thesis.
\end{proof}

\begin{lemma}\label{Iarrobino}
     Let $F_1,\ldots, F_r\in K[y_0, \ldots, y_s]$ be a set of $r$ {\it general} forms of fixed degree $d$. Then, for any positive integer number $i\leq d$, 
     \begin{equation*}
         \langle \frac{\partial^i F_1}{\partial y_0^{i_0}\dots \partial y_s^{i_s}},\ldots, \frac{\partial^i F_r}{\partial y_0^{i_0}\dots \partial y_s^{i_s}} \mid i_0+\dots i_s=i, i_j\geq 0 \rangle
     \end{equation*}
     is a $K$-vector space of dimension
     \begin{equation}
         \min\Set{r \binom{i+s}{s}, \binom{d+s-i}{s}}.
     \end{equation}
\end{lemma}
\begin{proof}
   It follows from  \cite[Proposition 3.4]{I}. 
\end{proof}
\medskip

 We explicitly note that Lemma \ref{Iarrobino} implies that $h_{max}$ is term-wise maximal. We also note that Proposition \ref{partial} extends Proposition 3.7 in \cite{FMMR22} and Theorem 3.5 in \cite{MRP} which refer to the case $m=2$.

We observe that $\alpha_0+\gamma_0=1$ and $\alpha_0+\beta_0=1+\binom{d+m-1}{m-1}$ so $h_0=1$, as expected. Moreover $\alpha_1+\gamma_1=m+n+1$ and $\alpha_1+\beta_1=m+\binom{d+m-2}{m-1}$, so $h_1=\alpha_1+\gamma_1$, and $\beta_1=\gamma_1$ if and only if $n+1=\binom{d+m-1}{m-1}$ which means that $A_F$ is a full Perazzo algebra.
\medskip

From now on we will  use the notation $s:=[\frac{d}{2}]$ so that $d=2s$ if it is even, and $d=2s+1$ if it is odd.  

We want to study  the unimodality of $h_{max}$. 
Note that in the range  $0\leq i\leq s$, $\alpha_i, \gamma_i$ are strictly increasing functions of $i$, independent of $d$, while $\beta_i$ is strictly decreasing and depends on $d$. 
\begin{lemma}\label{increasing}
    For any $i$,  $0\le i \le  s$, $\alpha_i+\gamma_i$ is a strictly increasing function of $i$,  while $\alpha_i+\beta_i$ is a strictly decreasing function of $i$.
\end{lemma}
\begin{proof}
    The first assertion is clear because both $\alpha_i$ and $\gamma_i$ are strictly increasing. To prove the second one, let $i\leq s-1$. We have:
    \begin{equation} 
    \begin{array}{ll} 
   &(\alpha_i+\beta_i)-(\alpha_{i+1}+\beta_{i+1})= \\[2ex]
   =&(\beta_i-\beta_{i+1})-(\alpha_{i+1}-\alpha_i) = \\[2ex]
   =&\binom{d+m-i-1}{m-1} -\binom{d+m-i-2}{m-1} -(\binom{m+i}{m-1}-\binom{m+i-1}{m-1})=\\[2ex]
   =&\binom{d+m-i-2}{m-2}-\binom{m+i-1}{m-2}>0.\\
    \end{array}\end{equation}
 Indeed from the hypothesis $i\leq s-1$ it follows $d-2i-1>0$, that is equivalent to $d+m-i-2>m+i-1$.   
\end{proof}

The maximal Hilbert vector of a Perazzo algebra with $m=2$, $n\ge 2$ and $d\ge 3$ is unimodal (see \cite[Theorem 3.6]{A} for the case $n=2$ and \cite[Theorem 4.12]{MRP} for the case $n\ge 2$). The result is no longer true for $m>2$ and in next theorem we will determine  when the maximal Hilbert vector of a Perazzo algebra with fixed
    $m, n, d$ and $m\geq 3$ is unimodal.

\begin{theorem}\label{non unimodal}
    Let $h_{max}$ be the maximal Hilbert vector of a Perazzo algebra with fixed
    $m, n, d$, $m\geq 3$. Let $s=[\frac{d}{2}].$  Then $h_{max}$ is unimodal if and only 
    \begin{enumerate}
        \item $\gamma_{s-1}<\beta_{s-1}$, and
        \item  $\alpha_{s-1}+\gamma_{s-1}\leq \alpha_s+\beta_s$.
    \end{enumerate}
\end{theorem}
\begin{proof}
     We will use repeatedly Lemma \ref{increasing}. 
     Assume first that conditions (1) and (2) are satisfied. From (1) it follows that $h_{s-1}=\alpha_{s-1}+\gamma_{s-1}$. We consider now $h_{s}$: if $h_s=\alpha_s+\gamma_s$, then $h_{max}$ is unimodal by Lemma \ref{increasing}; if $h_s=\alpha_s+\beta_s$, the unimodality of $h_{max}$ follows from Lemma \ref{increasing} and condition (2).

     Assume now that $h_{max}$ is unimodal. We observe that, by  Lemma \ref{increasing}, the existence of $i\leq s-1$ such that $\beta_i\leq \gamma_i$ is equivalent to $\beta_{s-1}\leq \gamma_{s-1}$. Therefore if, by contradiction $\beta_{s-1}\leq \gamma_{s-1}$, then there exists $i\leq s-1$ such that $h_i=\alpha_i+\beta_i$, i.e. $\beta_i\leq \gamma_i$. Being $\alpha_i+\beta_i\leq \alpha_i+\gamma_i$, using Lemma \ref{increasing} we get
     \begin{equation*}
         \alpha_{i+1}+\beta_{i+1}<\alpha_i+\beta_i\leq \alpha_i+\gamma_i < \alpha_{i+1}+\gamma_{i+1}.
     \end{equation*}
     Therefore $h_{i+1}=\alpha_{i+1}+\beta_{i+1}$ and $h_{i+1}<h_i=\alpha_i+\beta_i$ with $i+1\leq s$: this contradicts the unimodality of $h_{max}$. So condition (1) is satisfied. It implies that $h_{s-1}= \alpha_{s-1}+\gamma_{s-1}$. Finally, from $h_{s-1}\leq h_s$ it is immediate to deduce  condition (2).

\end{proof}

Now we want to translate the conditions (1) and (2) of Theorem \ref{non unimodal} in inequalities involving $n, m, d$, with $m\geq 3$.
We have to discute separately the cases $d$ even and $d$ odd.
\medskip

{\bf Condition (1), $d=2s$ even.}

$\gamma_{s-1}< \beta_{s-1}$ is equivalent to
\begin{equation*}
 \binom{s+m}{m-1}> (n+1)\binom{s+m-3}{m-1}.   
\end{equation*} This reduces to the inequality
\begin{equation*}
 ns^3-3(m-1)s^2-3[(m-1)^2+n]s-m(m-1)(m-2)< 0.   
\end{equation*}
Looking at the signs of the coefficients of the powers of $s$, from Descartes' rule of signs we deduce that the associated equation of degree $3$ in the unknown $s$ has at most one real positive solution $\bar s_1$. Therefore Condition (1) is never satisfied for $s$ large enough.
\medskip

{\bf Condition (1), $d=2s+1$ odd.}

$\gamma_{s-1}< \beta_{s-1}$ is equivalent to
\begin{equation*}
 \binom{s+m+1}{m-1}> (n+1)\binom{s+m-3}{m-1}.   
\end{equation*} This reduces to the inequality of degree $4$
\begin{equation*}
 ns^4+2(n-2m+2)s^3-(n+1+6m^2-6m-1)s^2-
 \end{equation*}
 \begin{equation*}
 (2n+2+4m^3-6m^2-2m+2)s-(m+1)m(m-1)(m-2)< 0.   
\end{equation*}
Again from Descartes' rule of signs we get that the associated equation has at most one positive solution $\bar s_2$, and we conclude as in the even case that Condition (1) is never satisfied for $s$ large enough.
\medskip

{\bf Condition (2), $d=2s$ even.}
The condition  $\alpha_{s-1}+\gamma_{s-1}\leq\alpha_s+\beta_s$ translates in  an inequality of degree $2$ in $s$, of the form
\begin{equation*}
ns^2-(3m+n-3)s-2(m-1)(m-2)\leq 0.
\end{equation*}
We conclude as in the previous cases.
\medskip

{\bf Condition (2), $d=2s+1$ odd.}
This time we get an inequality of degree $3$ in $s$ with at most one positive solution:
\begin{equation*}
    ns^3-4(m-1)s^2-(4m^2-8m+4+n)s-(m+1)(m-1)(m-2)\leq 0
\end{equation*}
and we conclude as in the previous cases.

Note that, even if we do not get explicit bounds on $s$,  we can summarize our computations in the following Theorem.
\begin{theorem}\label{non unimodal2}
   The maximal Hilbert vector of a Perazzo algebra with fixed  $n\geq m\geq 3$ is not unimodal for $d$ large enough.  
\end{theorem}
\begin{proof}
  It follows from the discussion after Theorem \ref{non unimodal}.  
\end{proof}
\medskip
\begin{example}
    
Case (2) in Theorem \ref{non unimodal} can fail as the following examples show. For $d$ even we take: $n=7$, $m=4$, $d=6$, $s=3$, $h_{max}=(1, 12, 42, 40, 42, 12, 1)$; $h_2=\alpha_2+\gamma_2$, $\alpha_2+\beta_2=10+\binom{9}{3}$, $h_3=\alpha_3+\beta_3$, $\alpha_3+\gamma_3=100$. An example with $d$ odd is the following: $n=13$, $m=3$, $d=5$, $s=2$, $h_{max}=(1, 17, 16, 17,1)$; $h_2=\alpha_2+\beta_2$ and $\alpha_2+\gamma_2=48$. 
\end{example}

\medskip

We want to study now if  Perazzo algebras with maximal $h$-vector have the WLP. In the special case $m=n=2$ this problem has been  solved in the negative in \cite{FMMR22}, and in the case $m=2$ and $n\ge 2$ in \cite{MRP}. 

As recalled in Section \ref{subsection Lefschetz}, it is well known that the $h$-vector of an Artinian graded algebra with the WLP is unimodal. We want to prove that for arbitrary $n,m$, even if the $h$-vector is unimodal, the Perazzo algebras with maximal $h$-vector fail the Weak Lefschetz Property.

\begin{theorem} \label{mainthm0}
Let $A_F$ be a Perazzo algebra with maximal $h$-vector for fixed $m, n, d$ with $n\geq m\geq 2$  and $d\geq 6$. Then $A_F$ fails the WLP.
\end{theorem}    
\begin{proof} For the case $m=2$ the reader can look at \cite{FMMR22} and \cite{MRP}. Let $m\geq 3$.
Let $h=(h_0, \ldots, h_d)$ be the $h$-vector of $A_F$. If $h$ is not unimodal, then the thesis trivially follows. So we assume that $h$ is unimodal.

Take $\ell$ a general linear form and consider the exact sequence appearing in Proposition \ref{exactsequence}:
\begin{equation}\label{seq}
0 \to (A_{\ell\circ F})(-1) \to A_F \to A_F/(\ell) \to 0.
\end{equation}
Since $h$ is unimodal we know from Theorem \ref{non unimodal} that
$\gamma_{s-1}< \beta_{s-1}^d$ and $\alpha_s+\beta_s^d\ge \alpha_{s-1}+\gamma_{s-1}$. Since $\beta_i$ depends on $d$, we keep track of it using the above notation.

We now discuss separately the cases $d$ odd and $d$ even.
\medskip

Assume first that $d=2s+1$ is odd. We observe that, from the assumption $n\geq m$, it follows $\beta_s^d< \gamma_s$; indeed $\beta_s^d=\binom{m+s}{m-1}$ and $\gamma_s=(n+1)\binom{m+s-2}{m-1}$, so a simple computation shows that $\beta_s^d< \gamma_s$ is equivalent to $s>\frac{m-1}{n}$.

Therefore the $h$-vector of $A_F$ is
\[ 
(\alpha_0+\gamma_0, \alpha_1+\gamma_1, \ldots, \alpha_{s-1}+\gamma_{s-1}, \alpha_s+\beta_s^d, \alpha_s+\beta_s^d, \alpha_{s-1}+\gamma_{s-1}, \ldots,  \alpha_0+\gamma_0).
\]

If $A_F$ has the WLP then the $h$-vector of $A_{\ell\circ F}$  will be 
\[
(\alpha_0+\gamma_0, \alpha_1+\gamma_1, \ldots, \alpha_{s-1}+\gamma_{s-1}, \alpha_s+\beta_s^d, \alpha_{s-1}+\gamma_{s-1}, \ldots, \alpha_1+\gamma_1, \alpha_0+\gamma_0).
\]
But $\alpha_s+\beta_s^d>\alpha_s+\beta_s^{d-1}$.  This is a contradiction because $\ell\circ F$ is a Perazzo polynomial of degree $d-1$ by Lemma \ref{lemmapartials2}: the condition $n+1\leq \binom{m+d-3}{m-1}$ is satisfied, otherwise a simple computation shows that $\beta_{s-1}^d\leq \gamma_{s-1}$, which implies that $h$ is not unimodal in view of Theorem \ref{non unimodal}.
\medskip

Let now $d=2s$ even. The $h$-vector of $A_F$ is
\[ 
(\alpha_0+\gamma_0, \alpha_1+\gamma_1, \ldots, \alpha_{s-1}+\gamma_{s-1}, h_s, \alpha_{s-1}+\gamma_{s-1}, \ldots,  \alpha_0+\gamma_0).
\]

If $A_F$ has the WLP then the $h$-vector of $A_{\ell\circ F}$  will be 
\[
(\alpha_0+\gamma_0, \alpha_1+\gamma_1, \ldots, \alpha_{s-1}+\gamma_{s-1}, \alpha_{s-1}+\gamma_{s-1}, \ldots, \alpha_1+\gamma_1, \alpha_0+\gamma_0),
\]
which implies that $\beta_{s-1}^{d-1}=\gamma_{s-1}$. But $\ell\circ F$ has odd degree $2s-1$ and is a Perazzo polynomial as in the previous case, hence $\beta_{s-1}^{d-1}<\gamma_{s-1}$: a contradiction. 
\end{proof}

\section{Minimal Hilbert function of a Perazzo algebra}\label{section minimal}
In this section we will compute the minimal $h$-vector of all Perazzo algebras with fixed $m,n,d$, extending the results obtained in \cite{FMMR22} and \cite{MRP}  for $m=2$. We will use Proposition \ref{rank}
 saying that, in the $h$-vector of the Perazzo algebra $A_F$, $h_i$ is equal to the rank of the matrix containing in the columns the coefficients of the partial derivatives of $F$ of order $i$.  So first of all we will give a precise description of this matrix. 

 \begin{proposition}\label{catalecticant}
     The matrix of the linear system defined in (\ref{conditions}) to compute $h_i$ has the following form
     \begin{equation}\label{matrix}
\begin{pmatrix}
0 & | & N_i   \\
 ---   & | & ---   \\
M_{i-1} &| & \Gamma_i  \\
\end{pmatrix}
\end{equation}
where:
\[
M_{i-1}= (C^0_{i-1} \ C^1_{i-1} \dots C^n_{i-1}),
\]
\[
N_i=\begin{pmatrix} C^0_i \\ 
C^1_i \\
\vdots \\
C^n_i,
\end{pmatrix},
\]
$C^k_i,  \Gamma_i$ are the catalecticant matrices for $p_k$, $k=0, \ldots, n$ and  $G$ defined as follows (see \cite{IK}, Definition 1.3):
\[
C^k_i=(p^k_{\delta+\eta})_{|\delta|=i, |\eta|=d-1-i} \text{ and } 
\ \Gamma_i=(G_{\delta+\eta})_{|\delta|=i, |\eta|=d-i}.
\]     
 \end{proposition}
 \begin{proof}
     It follows from the expression of $F$ and the assumption that $K$ has characteristic zero.
 \end{proof}

 We keep using the following notation introduced in Section \ref{maximal}:
 \begin{equation*}
 \alpha_i=\binom{m+i-1}{m-1}, \ \beta_i=\binom{d+m-i-1}{m-1}.    
\end{equation*}

\begin{theorem}\label{minimal}
    Let $m,n,d$ be fixed with $n\geq m\geq 2$.  Then the minimal $h$-vector of the  Perazzo algebras $A_F$, with $F$ polynomial of degree $d$ as in Definition \ref{Perazzo-form},  is $h_{\min}=(h_0, \ldots, h_d)$, where for $1\leq i\leq \frac{d}{2}$  
    \begin{equation}\label{minimum}
    h_i=
    \min\{2(n+1),\\
    \alpha_i+n+1, 
   \alpha_i+\beta_i \}.
    \end{equation}
\end{theorem}
\begin{proof}
   In view of Propositions \ref{rank} and \ref{catalecticant}, we have to look  for the Perazzo polynomials $F$ such that the rank of the matrix (\ref{matrix})  is minimal for any index $i=1, \ldots,[\frac{d}{2}]$. Therefore we can assume $G=0$, so that, for any $i$, $h_i=\rank M_{i-1}+\rank N_i$.   
   
   The minimal possible rank of each catalecticant matrix $C^k_{i-1}$ or $C^k_i$ is $1$. Therefore  the minimum between the number $n+1$ of catalecticant blocks of $M_{i-1}$ and $\beta_i=\binom{m+d-i-1}{m-1}=\dim K[U_1, \ldots, U_m]_{d-i}$ is a lower bound for the rank of $M_{i-1}$ for any $i\geq 1$. 
   
   Similarly a lower bound for the rank of $N_i$
is the minimum between the number of its columns, that is $\binom{m+i-1}{m-1}=\alpha_i$, and $n+1$ that is the number of its catalecticant blocks.

So we get the following lower bound: 
\[h_i\geq \min\Set{n+1, \beta_i}+\min\Set{n+1, \alpha_i}.\] Since $i\leq d-i$ in our range $i\leq \frac{d}{2}$, we get
\[h_i\geq \min\Set{2(n+1), n+1+\alpha_i, \alpha_i+\beta_i }.\]

To conclude the proof, we exhibit  an example of a Perazzo algebra with $h$-vector as in (\ref{minimum}). We observe that all the catalecticant matrices of the polynomial  $L^d$, where $L$ is a linear form, have rank 1. Therefore, in view of  \cite[Corollary 3.2]{I} it is enough to take $F=X_0L_0^{d-1}+\dots + X_nL_n^{d-1}$, where $L_0, \ldots, L_n\in K[U_1, \ldots, U_m]_1$ are general linear forms.

\end{proof}
For $m=2$ formula (\ref{minimum}) gives the expression of the minimum $h$-vector found in  \cite[Proposition 2.8]{FMMR22} and \cite[Theorem 3.4]{MRP}.

\medskip

 In \cite{FMMR22}, \cite{A} and \cite{MRP} it was proved that in the case $m=2$  the minimal $h$-vector of Perazzo algebras is always unimodal and that the Perazzo algebras with minimal $h$-vector have the Weak Lefschetz Property  provided that $d\geq 2n$. This does not always happen in the general case, as next example shows that $h_{min}$ can be non unimodal for some integers $n, m, d$.

\begin{example}
    Let $m=3$, $n=9$, $d=4$. These are the invariants of  the famous Stanley's example (\cite{S}),  whose  $h$-vector is $(1, 13, 12,13,1)$: it is clearly non unimodal and it is minimal for these invariants (see \cite{CGIZ}); it corresponds to a full Perazzo algebra. 
\end{example}

\begin{theorem} \label{unimodality2}
    The minimal $h$-vector of  Perazzo algebras with invariants $n\geq m> 2$  is unimodal if and only if  $n+1\leq \beta_{s-1}$,  where $s=[\frac{d}{2}]$.
    
\end{theorem}
\begin{proof}  
    Theorem \ref{minimal} implies that, given $n, m, d$, for any $1\leq i\leq s$ with $s=[\frac{d}{2}]$:
\begin{enumerate}
    \item if $n+1\leq \alpha_i\leq \beta_i$, then $h_i=2(n+1)$ ;
    \medskip
    
    \item if \ $\alpha_i< n+1 \leq \beta_i$, then $h_i=n+1+\alpha_i$, which is an \emph{increasing} function of $i$;
    \medskip
    
    \item if $n+1>\beta_i$, then $h_i=\alpha_i+\beta_i$, which is a \emph{decreasing} function of $i$ depending also on $d$. 
\end{enumerate}
It follows that $h_{min}$ is unimodal if and only if $h_{s-1}$ is not of the form (3). This proves the Theorem.
\end{proof}

We are now able to prove that  the Weak Lefschetz Property holds for the Perazzo algebras with minimal $h$-vector, provided it  satisfies the  condition that $n+1\leq \beta_s$, which implies unimodality.

\begin{theorem}\label{minimal with WLP}
    Let $A_F$ be a Perazzo algebra  with minimal $h$-vector for fixed
$m, n, d$ with $n \geq m \geq 3$. Assume that the $h$-vector of $A_F$ is unimodal and that $n+1\leq \beta_s$, where $s=[\frac{d}{2}]$. Then $A_F$ has the WLP.
\end{theorem}
\begin{proof}
    To prove that $A_F$ has the WLP it is enough to check that, for a general linear form $\ell$, the multiplication map  $\times \ell: (A_F)_{s-1} \to (A_F)_s$ is injective. If by contradiction it is not injective, then using Proposition \ref{exactsequence} and Lemma \ref{lemmapartials2}, we get that $\dim (A_{\ell\circ F})_{s-1}<\dim (A_F)_{s-1}$. But our assumption on $n,m,d$ implies that the component of index $s-1$ of the minimal $h$-vector is the same for degrees $d$ and $d-1$. This contradicts the minimality of the $h$-vector of $A_F$.
\end{proof}
Note that if $d$ is odd, then the condition $n+1\leq \beta_s$ means that the $h$-vector of $A_{\ell\circ F}$ is unimodal.

We characterize now the integers $n, m, d$ such that the Perazzo algebras with these invariants have all the same  Hilbert function, i.e. the maximal and the minimal $h$-vector coincide.

\begin{proposition}
    Let $n,m,d$ be positive integers with $n\geq m\geq 2$, $n+1 \leq \binom{d+m-2}{m-1}$. Then the Perazzo algebras with these invariants have $h_{min}=h_{max}$ if and only if 
    \begin{equation}\label{max=min}
        \binom{d+m-3}{m-1}\leq n+1.
    \end{equation}
\end{proposition}
\begin{proof}
    We will prove that Equation \ref{max=min} is equivalent to $h_{i,min}=h_{i,max}=\alpha_i+\beta_i$ for any index $i$.
   We will then exclude the possibility that $h_{2,min}=h_{2,max}\neq \alpha_2+\beta_2.$

   Since $\alpha_i+\beta_i$ is a decreasing function of $i$, the first assertion is equivalent to $h_{2,min}=h_{2,max}=\alpha_2+\beta_2.$ But $h_{2,min}=\alpha_2+\beta_2$ if and only if $\binom{m+d-3}{m-1}\leq n+1$ and $h_{2,max}=\alpha_2+\beta_2$ if and only if $\binom{m+d-3}{m-1}\leq (n+1)m$, which proves our claim.

   We assume now by contradiction that $h_{2,min}=h_{2,max}\neq \alpha_2+\beta_2.$ This means that $h_{2,max}=\alpha_2+\gamma_2<\alpha_2+\beta_2$, and $h_{2,min}=\min\{ 
   2n+2, n+1+\binom{m+1}{2}\}<\alpha_2+\beta_2.$
   If $n+1\leq \binom{m+1}{2}$, then $2n+2=\binom{m+1}{2}+(n+1)m$ which is equivalent to $(n+1)(2-m)=\binom{m+1}{2}$, but this  is impossible because $m\geq 2$. If $n+1\geq \binom{m+1}{2}$, then we would have $n+1=(n+1)m$: contradiction.
   \end{proof}

\begin{remark}
Perazzo algebras with  $h_{min}=h_{max}$ clearly include full Perazzo algebras. It has been conjectured in \cite[Conjecture 2.6]{BGIZ} that the $h$-vectors of full Perazzo algebras are minimal among the $h$-vectors of all Artinian Gorenstein algebras with the same degree and codimension. This conjecture has been proved  for $d=4$ and $m=3,4,5$ in \cite{CGIZ}, and for any degree and $m=3$ in \cite{BGIZ}.  

\end{remark}
\section{Minimal free resolution  of a Perazzo algebra}
\label{MFRPerazzo}

In this section, we determine the minimal free resolution of an Artinian Gorenstein algebra corresponding to a Perazzo threefold of  $\PP^4$ with termwise minimal Hilbert function, i.\,e.\ of the following  type $(1,5,6,\dots,6,5,1)$. 

When we deal with  Perazzo hypersurfaces in $\PP^4$,  we use the notations 
$S=K[X,Y,Z,U,V]$ and $R=K[x,y,z,u,v]$. We have
\begin{equation}
\label{Perazzo-form}
F=Xp_0+Yp_1+Zp_2+G \quad \text{where} \ p_0, p_1, p_2, G  \in K[U,V]
\end{equation}
and any choice of $p_0, p_1, p_2$ will be algebraically dependent. 

An explicit classification of the possible dual generators $F$ of degree $d \ge 5$ defining a Perazzo threefold with termwise minimal Hilbert function is given in \cite[Theorem 5.4]{FMMR22}. We state a slightly rephrased version of this result. 

\begin{lemma}\label{lemma:min_perazzo}
    Let $F \in K[X,Y,Z,U,V]$ be a Perazzo form such that the algebra $A_F$ has minimal Hilbert function $(1,5,6, 6, \ldots, 6,5,1)$. Then the dual generator $F$ can be expressed as
    \begin{enumerate}[(i)]
    \item \label{minimalHFPerazzo1}$XU^{d-1}+YU^{d-2}V+ZU^{d-3}V^2$, 
    \item \label{minimalHFPerazzo2}$XU^{d-1}+YU^{d-2}V+ZV^{d-1}$, or
    \item \label{minimalHFPerazzo3} $XU^{d-1}+Y(U+\lambda V)^{d-1}+ZV^{d-1}$ with $\lambda\in K^*$
\end{enumerate}
after a linear change of variables. 
\end{lemma}
\begin{proof}
    By \cite[Theorem 5.4]{FMMR22}, there are three classes of forms up to a linear change of variables. In the first case, $F$ can be written as
    \begin{align*}
        F&=XU^{d-1}+YU^{d-2}V+ZU^{d-3}V^2+aU^d+bU^{d-1}V+cU^{d-2}V\\
        &=(X+aU)U^{d-1}+(Y+bU)U^{d-2}V+(Z+cU)U^{d-3}V^2
    \end{align*}
    for $a,b,c\in K$. A further linear change of variables gives the case (i). The other cases are obtained in the same manner.
\end{proof}
 
We know by \cite[Theorem 4.3]{FMMR22} that, in all the three cases of Lemma \ref{lemma:min_perazzo}, $A_F$ has the WLP which allows us to prove the following key lemma.

\begin{lemma}\label{tor=0}
 Let $A_F$ be  an Artinian Gorenstein algebra associated to a Perazzo threefold $X=V(F)$ of  $\PP^4$ of degree $d\ge 5$ and  with termwise minimal Hilbert function. Let $\ell \in A_F$
 be a general linear form and set $B=A_F/(\ell)$. We have:
 $$
 \beta _{ij}^R(B):=[Tor^R_{i}(B,K)]_j=0 \text{ for }
 j>i+2. $$
 \end{lemma}
\begin{proof}
 By \cite[Theorem 4.3]{FMMR22}  the algebra $A_F$ has the WLP which implies that  the $h$-vector of $B$ is $(1,4,1)$ and the socle degree of $B$ is 2, i.e. 
 $\reg(B)=2$.  Therefore, the minimal graded free resolution of $B$ as $R$-module has the following shape:
 $$
 0\rightarrow  \begin{array}{c}
R(-6)^{\beta_{56}^R(B)} \\
\oplus \\
R(-7)^{\beta_{57}^R(B)}
\end{array}\rightarrow \begin{array}{c}
R(-5)^{\beta_{45}^R(B)} \\
\oplus \\
R(-6)^{\beta_{46}^R(B)}
\end{array}\rightarrow \begin{array}{c}
R(-4)^{\beta_{34}^R(B)} \\
\oplus \\
R(-5)^{\beta_{35}^R(B)}
\end{array}\rightarrow $$
$$\begin{array}{c}
R(-3)^{\beta_{23}^R(B)} \\
\oplus \\
R(-4)^{\beta_{24}^R(B)}
\end{array}\rightarrow \begin{array}{c}
R(-1)\\ \oplus \\
R(-2)^9 \\
\oplus \\
R(-3)^{\beta_{13}^R(B)}
\end{array}\rightarrow R \rightarrow B \rightarrow 0
 $$
and we conclude that $
 \beta _{ij}^R:=[Tor^R_{i}(B,K)]_j=0$   for  $j>i+2 $ which proves what we want. 
\end{proof}

\begin{example}\label{ex_ind}
Using Macaulay2 \cite{M2},  we have computed the Betti table of Artinian Gorenstein algebras $A_F$ associated to  all 3 possible types of Perazzo threefolds $F$ of  $\PP^4$ with termwise minimal Hilbert function and degree $5\le d\le 8$ (see Lemma \ref{lemma:min_perazzo}),

{\scriptsize
\begin{verbatim}
				+-----------------------+
				|       0  1  2  3  4  5|
				|total: 1 14 35 35 14  1|
				|    0: 1  .  .  .  .  .|
				|    1: .  9 17 12  3  .|
				|    2: .  1  3  3  1  .|
				|    3: .  1  3  3  1  .|
				|    4: .  3 12 17  9  .|
				|    5: .  .  .  .  .  1|
				+-----------------------+
\end{verbatim}
\begin{verbatim}
				+-----------------------+
				|       0  1  2  3  4  5|
				|total: 1 14 35 35 14  1|
				|    0: 1  .  .  .  .  .|
				|    1: .  9 17 12  3  .|
				|    2: .  1  3  3  1  .|
|    3: .  .  .  .  .  .|
				|    4: .  1  3  3  1  .|
				|    5: .  3 12 17  9  .|
				|    6: .  .  .  .  .  1|
				+-----------------------+
\end{verbatim}
\begin{verbatim}
				+-----------------------+
				|       0  1  2  3  4  5|
				|total: 1 14 35 35 14  1|
				|    0: 1  .  .  .  .  .|
				|    1: .  9 17 12  3  .|
				|    2: .  1  3  3  1  .|
|    3: .  .  .  .  .  .|
|    4: .  .  .  .  .  .|
				|    5: .  1  3  3  1  .|
				|    6: .  3 12 17  9  .|
				|    7: .  .  .  .  .  1|
				+-----------------------+
\end{verbatim}

\begin{verbatim}
				+-----------------------+
				|       0  1  2  3  4  5|
				|total: 1 14 35 35 14  1|
				|    0: 1  .  .  .  .  .|
				|    1: .  9 17 12  3  .|
				|    2: .  1  3  3  1  .|
|    3: .  .  .  .  .  .|
|    4: .  .  .  .  .  .|
|    5: .  .  .  .  .  .|
				|    6: .  1  3  3  1  .|
				|    7: .  3 12 17  9  .|
				|    8: .  .  .  .  .  1|
				+-----------------------+
\end{verbatim}
}
    
\end{example}

\begin{theorem} \label{main}  Let $A_F$ be  an Artinian Gorenstein algebra associated to a Perazzo threefold $X=V(F)$ of  $\PP^4$ of degree $d\ge 5$ and  with minimal Hilbert function. The Betti diagram of $A_F$ looks like

\begin{verbatim}
				+-----------------------+
				|       0  1  2  3  4  5|
				|total: 1 14 35 35 14  1|
				|    0: 1  .  .  .  .  .|
				|    1: .  9 17 12  3  .|
				|    2: .  1  3  3  1  .|
|    3: .  .  .  .  .  .|
|    .:.   .  .  .  .  .|
|  d-3: .  .  .  .  .  .|
				|  d-2: .  1  3  3  1  .|
				|  d-1: .  3 12 17  9  .|
				|    d: .  .  .  .  .  1|
				+-----------------------+
\end{verbatim}
\end{theorem}
\begin{proof}
    We  proceed by induction on $d$. For $5\le d\le 8$ the result is true (see Example \ref{ex_ind}). Assume $\deg(F)=d+1\ge 9$. Let $\ell \in A_F$ be a general linear form and consider the exact sequence:
    $$ \label{exact}
    0 \longrightarrow A_{\ell\circ F}  \longrightarrow A_F  \longrightarrow B=A_F/(\ell)  \longrightarrow 0
    $$
    which gives us the long exact sequence
    $$
    \cdots \longrightarrow [Tor^R_{i+1}(B,K)]_j \longrightarrow [Tor^R_{i}(A_{\ell \circ F}(-1),K)]_j  \longrightarrow [Tor^R_{i}(A_F,K)]_j $$
    $$\longrightarrow [Tor^R_{i}(B,K)]_j \longrightarrow \cdots .
    $$
    Using Lemma \ref{tor=0}, we get that for any $j>i+3$ we have:
  \begin{equation}\label{keyiso} \begin{array}{rcl} \beta_{ij-1}^R(A_{\ell\circ F}) & = &  [Tor^R_{i}(A_{\ell \circ F},K)]_{j-1} \\
  & = & [Tor^R_{i}(A_{\ell \circ F}(-1),K)]_j  \\
  & = & [Tor^R_{i}(A_F,K)]_j \\
  & = & \beta_{ij}^R(A_F) .
  \end{array}\end{equation}
   By Lemma \ref{lemmapartials2}, $A_{\ell \circ F}$ is an Artinian  Gorenstein algebra associated to a Perazzo form $\ell \circ F$ of degree $d$. Using the exact sequence (\ref{exact})
 and the fact that $A_F$ has the WLP 
 (\cite[Theorem 4.3]{FMMR22})  we obtain that the Hilbert function of $A_{\ell \circ F}$ is  termwise minimal, i.e. $A_{\ell \circ F}$ has Hilbert function $(1,5,6,6,\cdots ,6,6,5,1)$. 
 Therefore, by hypothesis of induction we know all graded Betti numbers of $A_{\ell \circ F}$ and the equalities (\ref{keyiso}) give us:
$$ \begin{array}{cclcl}
1& = & \beta_{1,d-1}(A_{\ell \circ F}) & = & \beta _{1,d}(A_F) \\
3& = & \beta_{1,d}(A_{\ell \circ F}) & = & \beta _{1,d+1}(A_F) \\
0& = & \beta_{1,j-1}(A_{\ell \circ F}) & = & \beta _{1,j}(A_F)  \text{ for } d+j\ge 5;    
\end{array}
$$

$$ \begin{array}{rclcl}
3& = & \beta_{2,d}(A_{\ell \circ F}) & = & \beta _{2,d+1}(A_F) \\
12& = & \beta_{2,d+1}(A_{\ell \circ F}) & = & \beta _{2,d+2}(A_F) \\
0& = & \beta_{2,j-1}(A_{\ell \circ F}) & = & \beta _{2,j}(A_F)  \text{ for } d+j\ge 6;
\end{array}
$$

 $$ \begin{array}{rclcl}
3& = & \beta_{3,d+1}(A_{\ell \circ F}) & = & \beta _{3,d+2}(A_F) \\
17& = & \beta_{3,d+2}(A_{\ell \circ F}) & = & \beta _{3,d+3}(A_F) \\
0& = & \beta_{3,j-1}(A_{\ell \circ F}) & = & \beta _{3,j}(A_F)  \text{ for } d+j\ge 7; 
\end{array}
$$

$$ \begin{array}{rclcl}
1& = & \beta_{4,d+2}(A_{\ell \circ F}) & = & \beta _{4,d+3}(A_F) \\
9& = & \beta_{4,d+3}(A_{\ell \circ F}) & = & \beta _{4,d+4}(A_F) \\
0& = & \beta_{4,j-1}(A_{\ell \circ F}) & = & \beta _{4,j}(A_F)  \text{ for } d+j\ge 8.  
\end{array}
$$

Therefore, the Betti diagram of $A_F$ being $\deg(F)=d+1$ has the following shape ($*$ means not yet determined and $.$ means zero):
\begin{verbatim}
				+-----------------------+
				|       0  1  2  3  4  5|
				|total: 1 14 35 35 14  1|
				|    0: 1  .  .  .  .  .|
				|    1: .  *  *  *  *  .|
				|    2: .  *  *  *  *  .|
|    3: .  *  *  *  *  .|
|    4:.   *  *  *  *  .|
|    5:.   .  .  .  .  .|
|    .:.   .  .  .  .  .|
|  d-2: .  .  .  .  .  .|
				|  d-1: .  1  3  3  1  .|
				|    d: .  3 12 17  9  .|
				|  d+1: .  .  .  .  .  1|
				+-----------------------+
\end{verbatim}

Using now the fact  that the minimal graded free R-resolution of an Artinian Gorenstein algebra $A_F$ is self dual we conclude what we want.
\end{proof}

\section{Final remarks and open problems}\label{final}
We end this paper with a couple of concrete problems which naturally arise from our results and we believe they deserve further consideration.

In Theorem \ref{mainthm0} we prove that any Perazzo algebra $A_F$ with maximal $h$-vector for fixed $n\ge m\ge 2$ and $d\ge 6$ fails WLP while in Theorem \ref{non unimodal} we determine  all Perazzo algebras with unimodal maximal $h$-vector for fixed $n$, $m$ and $d$ with $m\ge 3$. The analogous results for  Perazzo algebras $A_F$ with minimal $h$-vector for fixed  $n$, $m$ and $d$ with $n\ge m\ge 3$ are obtained in Theorems \ref{minimal with WLP} and \ref{unimodality2}.

For Perazzo algebras $A_F$ with intermediate $h$-vector both possibilities occur: there are examples failing WLP and examples satisfying WLP as well as examples with unimodal $h$-vector and examples with non-unimodal $h$-vector (see \cite{FMMR22}).
Therefore,  the major questions/problems left open are the following two:

\begin{problem}\label{wlp-unimodality}
(i) To classify all Perazzo algebras $A_F$ with unimodal Hilbert function.

(ii) To classify all Perazzo algebras $A_F$ with WLP.
    \end{problem}

For a complete answer to Problem \ref{wlp-unimodality} for $n=m=2$ the reader can look at \cite{A} and \cite{FMMR22} and for $n\ge2$ and $m=2$ at \cite{MRP}. To our knowledge for all other values $n\ge m\ge 3$ no answer is known.

\vskip 4mm
In the last decades big effort has been made in understanding the minimal free resolution (MFR, for short) of any artinian Gorenstein algebra. In 1977,  Buchsbaum and Eisenbud proved that any Gorenstein codimension 3 ideal is generated by the $2t \times 2t$ pfaffians of a skew symmetric matrix of size $(2t+1)\times (2t+1)$ and this fact  completely determines the MFR of any Artinian Gorenstein algebra of codimension 3 (see \cite[Theorem 2.1]{BE}). Nevertheless for codimension $\ge 4$ little is known apart from the selfduality (up to twist) of the MFR of any Artinian Gorenstein algebra. Using our knowledge of Perazzo algebras  we propose as an intermediate step the following problem:

\begin{problem}\label{resolution}  

(i) To determine the MFR of any Perazzo algebra with minimal (resp. maximal) Hilbert function.

(ii) To determine the MFR of any full Perazzo algebra.

(iii) To determine the MFR of {\em any} Perazzo algebra.
\end{problem}

The above problem is interesting {\it per se} but also because we believe  that the WLP of Perazzo algebras of any codimension $c$ could be determined by their MFR.
\providecommand{\bysame}{\leavevmode\hbox to3em{\hrulefill}\thinspace}
\providecommand{\MR}{\relax\ifhmode\unskip\space\fi MR }
\providecommand{\MRhref}[2]{%
  \href{http://www.ams.org/mathscinet-getitem?mr=#1}{#2}
}

\end{document}